\newif\ifacm
\acmfalse
\acmtrue

\ifacm
\documentclass{sig-alternate-05-2015}
\usepackage{etoolbox}
\newtoggle{acm}
\toggletrue{acm}
\newtoggle{public}
\toggletrue{public}
\newtoggle{doubleblind}
\toggletrue{doubleblind}
\togglefalse{doubleblind}

\else

\documentclass[12pt]{article}
\usepackage[margin=2cm]{geometry}
\usepackage{etoolbox}
\newtoggle{acm}
\togglefalse{acm}
\newtoggle{public}
\toggletrue{public}
\newtoggle{doubleblind}
\togglefalse{doubleblind}

\fi

\iftoggle{acm}
{
\usepackage[                   
  backend=bibtex8,             
  bibencoding=utf8,            
  style=numeric,            
  natbib=true,                 
  hyperref=true,               
  backref=false,                
  isbn=false,                  
  url=false,                    
  doi=false,                    
  urldate=long,                
  firstinits=true,
  terseinits=false,
  maxnames=3,%
  minnames=1,%
  maxbibnames=25,%
  minbibnames=15,%
  maxcitenames=2,%
  mincitenames=1%
]{biblatex}
}
{
\usepackage[                   
  backend=bibtex8,             
  bibencoding=utf8,            
  style=alphabetic,            
  natbib=true,                 
  hyperref=true,               
  backref=true,                
  isbn=false,                  
  url=true,                    
  doi=true,                    
  urldate=long,                
  firstinits=true,
  terseinits=false,
  maxnames=3,%
  minnames=1,%
  maxbibnames=25,%
  minbibnames=15,%
  maxcitenames=2,%
  mincitenames=1%
]{biblatex}
}

\renewbibmacro*{doi+eprint+url}{%
  \iftoggle{bbx:doi}
    {\printfield{doi}}
    {}%
  \newunit\newblock
  \iftoggle{bbx:eprint}
    {\iffieldundef{doi}{\iffieldundef{url}{\printfield{eprint}}{}}{}}
    {}%
  \newunit\newblock
  \iftoggle{bbx:url}
    {\iffieldundef{doi}{\printfield{url}}{}}
    {}}

\usepackage[
    hypertexnames,%
    citecolor=blue,%
    colorlinks=true,%
    linkcolor=red%
]{hyperref}
\usepackage[all]{hypcap}

\usepackage[font=scriptsize,caption=false]{subfig}
\iftoggle{acm}{
}
{}

\iftoggle{public}{
\usepackage[switch]{lineno}
}
{}
\usepackage{amssymb,amsmath}
\usepackage{dsfont} 
\usepackage[mathscr]{eucal}

\usepackage{verbatim}
\usepackage{color}
\usepackage{graphicx}

\usepackage{tikz}
\usetikzlibrary{arrows}
\usetikzlibrary{calc}

\usepackage{xcolor}
\definecolor{darkblue}{HTML}{0000AD}
\definecolor{darkgreen}{HTML}{008600}
\definecolor{darkred}{HTML}{8B0000}
\definecolor{darkgray}{HTML}{666666}
\definecolor{_mage}{HTML}{912830}
\definecolor{_cyan}{HTML}{31837a}
\definecolor{_purp}{HTML}{49425c}

\usepackage{url}

\usepackage{environ}
\NewEnviron{killcontents}{}

\newcommand{\Poincare}{Poincar\'{e} }

\newcommand{\figref}[1]{{Fig.~\ref{#1}}}
\newcommand{\fig}[1]{{\figref{fig:#1}}}
\newcommand{\secref}[1]{{Sec.~\ref{#1}}}
\newcommand{\sct}[1]{{\secref{sec:#1}}}

\newcommand{\see}[1]{(see~#1)}

\newcommand{\set}[1]{\left\{ #1 \right\}}
\newcommand{\paren}[1]{\left( #1 \right)}
\newcommand{\brak}[1]{\left[ #1 \right]}
\newcommand{\abs}[1]{\left| #1 \right|}
\newcommand{\ip}[2]{\langle #1,#2 \rangle}
\newcommand{\ipM}[2]{\ip{#1}{#2}_M}

\newcommand{\norm}[1]{\left\| #1 \right\|}

\newcommand{\mat}[2]{\brak{\begin{array}{#1} #2 \end{array}}}

\newcommand{\td}[1]{\widetilde{#1}}

\newcommand{\eqnn}[1]{\begin{equation}\begin{aligned} #1 \end{aligned}\end{equation}}

\newcommand{\sm}{\setminus}
\newcommand{\into}{\rightarrow}
\newcommand{\goesto}{\rightarrow}

\newcommand{\tr}{\top}

\newcommand{\I}{I_d}
\newcommand{\II}{I_{2d}}

\newcommand{\R}{\mathbb{R}}

\newcommand{\N}{\mathbb{N}}

\newcommand{\e}{\mathscr}
\newcommand{\q}{p}
\newcommand{\p}{q}
\newcommand{\w}{\omega}

\newcommand{\vphi}{\varphi}
\newcommand{\veps}{\varepsilon}

\newtheorem{proposition}{Proposition}
\newtheorem{definition}{Definition}
\newtheorem{theorem}{Theorem}
\newtheorem{corollary}{Corollary}
\newtheorem{lemma}{Lemma}
\newtheorem{claim}{Claim}
\newtheorem{assumption}{Assumption}
\newtheorem{remark}{Remark}
\newtheorem{example}{Example}

\newcommand{\defna}[2]{\begin{definition}[#1] #2 \end{definition}}
\newcommand{\rem}[1]{\begin{remark} #1 \end{remark}}
\newcommand{\remna}[2]{\begin{remark}[#1] #2 \end{remark}}

\newcommand{\asmpna}[2]{\begin{assumption}[#1] #2 \end{assumption}}

\newcommand{\lemna}[2]{\begin{lemma}[#1] #2 \end{lemma}}

\newcommand{\thmna}[2]{\begin{theorem}[#1] #2 \end{theorem}}

\newcommand{\asmpdiff}{Assump.~\ref{asmp:diff} ($C^r$ vector field and reset map)}

\newcommand{\asmpindep}{Assump.~\ref{asmp:indep} (independent constraints)}

\newcommand{\asmporth}{Assump.~\ref{asmp:orth} (orthogonal constraints)}

\newcommand{\asmpflow}{Assump.~\ref{asmp:flow} (existence and uniqueness of flow)}

\newcommand{\remindep}{Remark~\ref{rem:indep} (independent constraints)}

\newcommand{\remtrj}{Remark~\ref{rem:trj} (admissible trajectories)}
\newcommand{\remadact}{Remark~\ref{rem:adact}}
\newcommand{\defcontact}{Def.~\ref{def:contact} (contact modes)}

\newcommand{\defadact}{Def.~\ref{def:adact} (admissible constraint activation/deactivation)}
\newcommand{\deftrj}{Def.~\ref{def:trj} (admissible trajectory)}
\newcommand{\defseq}{Def.~\ref{def:seq} (contact mode sequence)}

\newcommand{\lemcont}{Lem.~\ref{lem:cont} (continuity across contact mode sequences)}
\newcommand{\thmdiffb}{Thm.~\ref{thm:diffb} (piecewise differentiability across contact mode sequences)}

\newcommand{\Frechet}{Fr{\'{e}}chet}

\title{Piecewise--differentiable trajectory outcomes in mechanical systems subject to unilateral constraints}

\iftoggle{doubleblind}
{ 
	\author{Omitted for submission}
}
{

\author{
Andrew~M.~Pace$^*$%
\and Samuel~A.~Burden%
\thanks{Department of Electrical Engineering, University of Washington, Seattle, WA, USA ({\tt apace2,sburden@uw.edu}).}%
\thanks{This work was supported by ARO Young Investigator Program Award \#W911NF-16-1-0158 to S. Burden.}
}
}

\date{}

\begin{document}

\maketitle

\iftoggle{public}
{
}
{}

\abstract{

We provide conditions under which trajectory outcomes in mechanical systems subject to unilateral constraints depend piecewise--differentiably on initial conditions, even as the sequence of constraint activations and deactivations varies.
This builds on prior work that provided conditions ensuring existence, uniqueness, and continuity of trajectory outcomes, and extends previous differentiability results that applied only to fixed constraint (de)activation sequences.
We discuss extensions of our result and implications for assessing stability and controllability.
}

\iftoggle{acm}
{
}

\section{Introduction}
\label{sec:intro}

\begin{figure}[h!]
\centering
{
\iftoggle{acm}
{
\hfill%
{
\def\svgwidth{0.5\columnwidth} 
\resizebox{.45\columnwidth}{!}{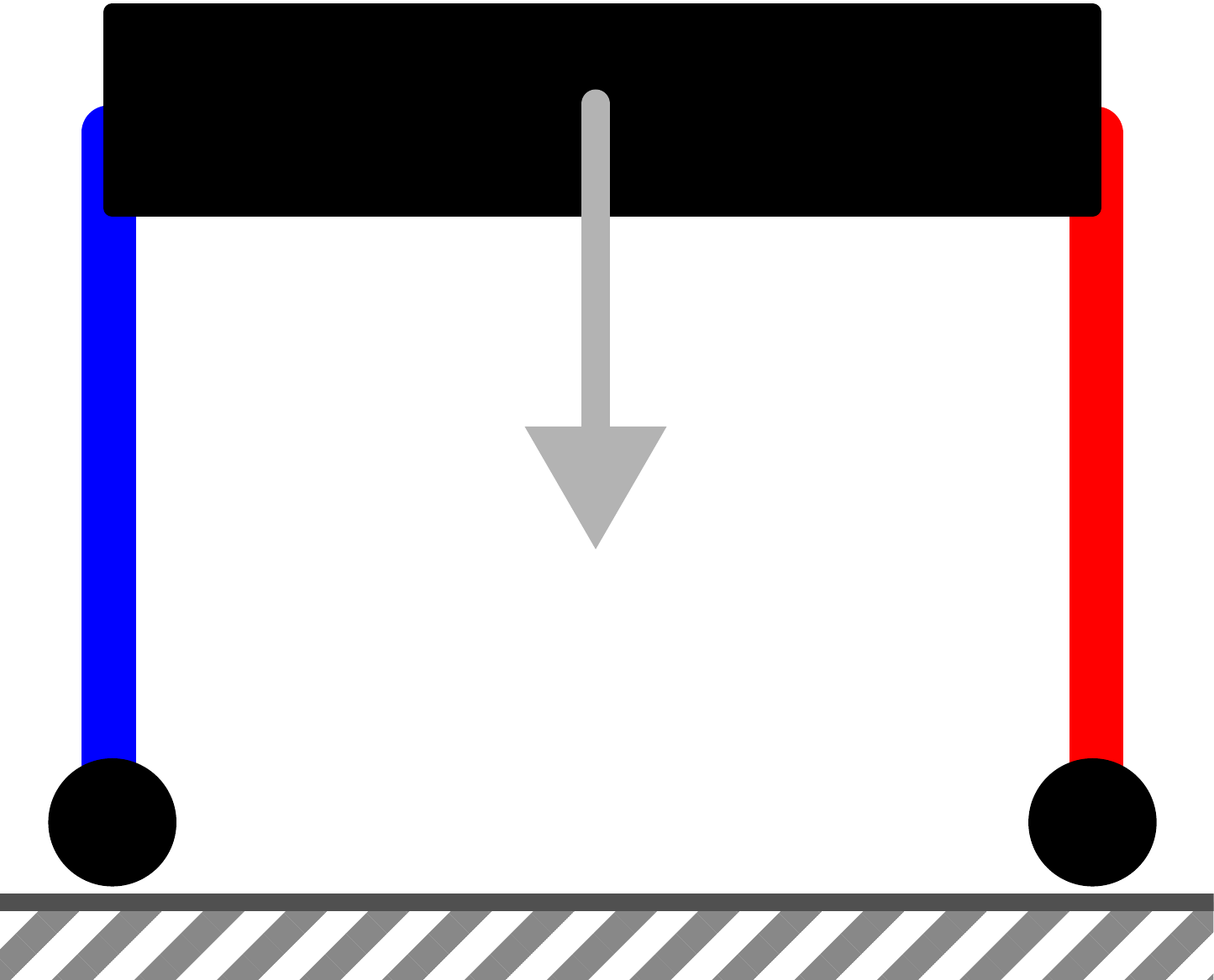}
}
\hfill%
{
\def\svgwidth{0.5\columnwidth} 
\resizebox{.45\columnwidth}{!}{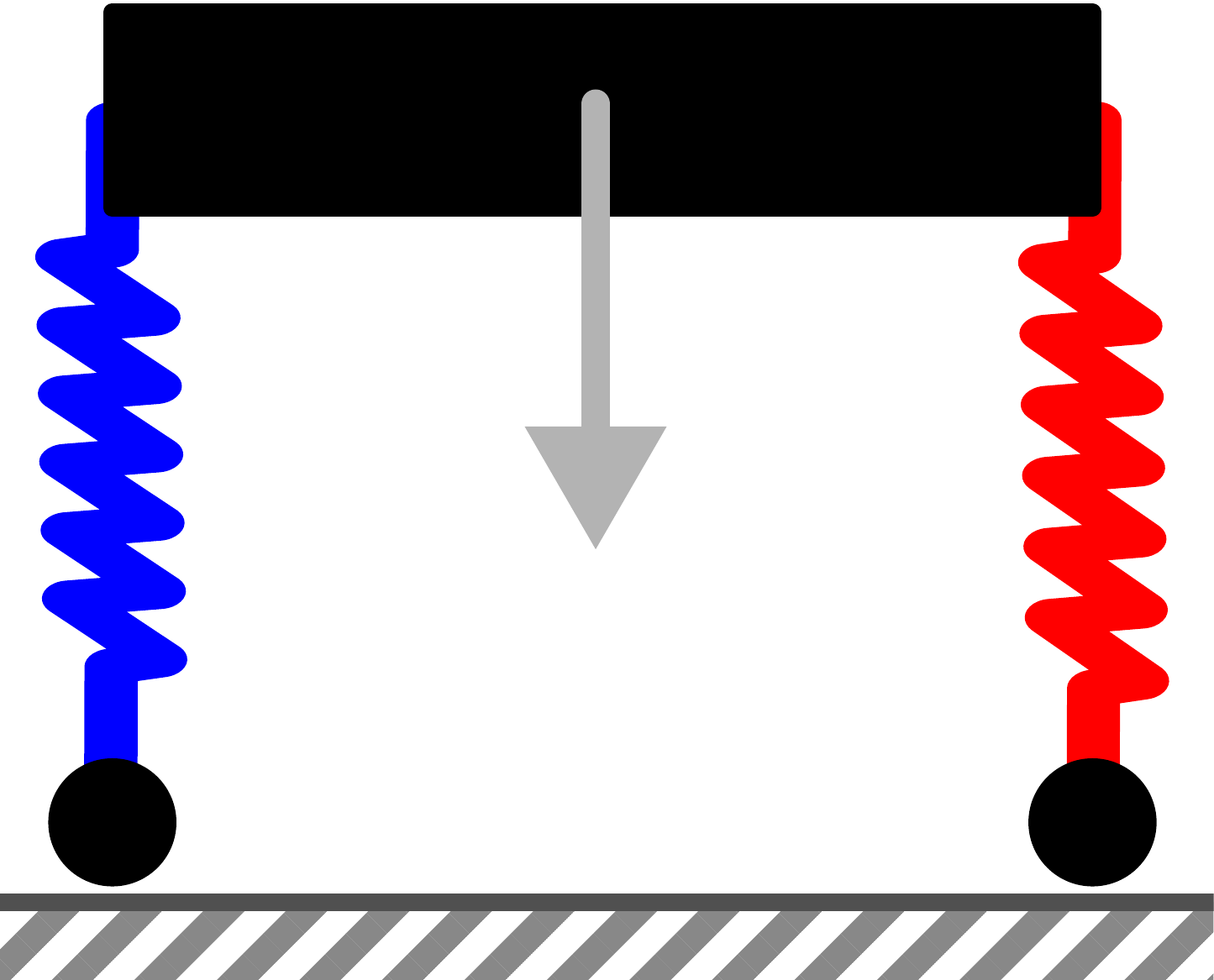}
}
\includegraphics[width=\columnwidth]{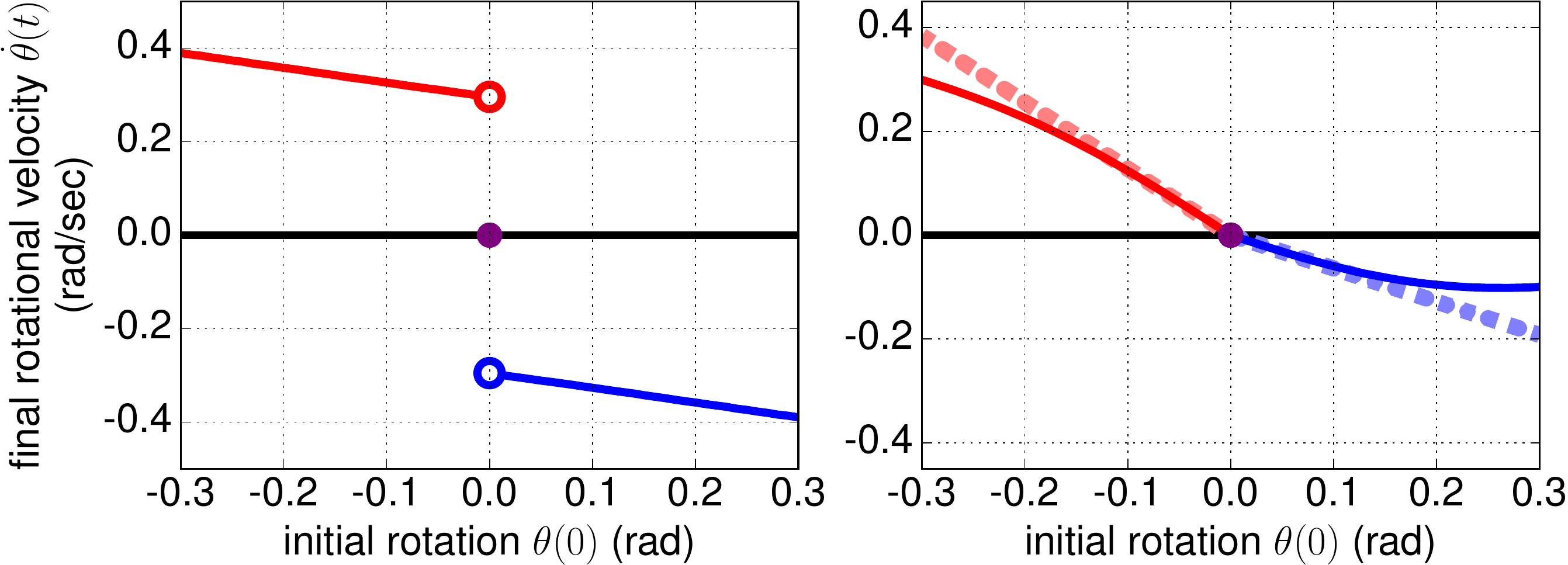}
}
{
\quad\includegraphics[width=.4\columnwidth]{fig/trot_td_stiff.pdf}\quad\quad
\includegraphics[width=.4\columnwidth]{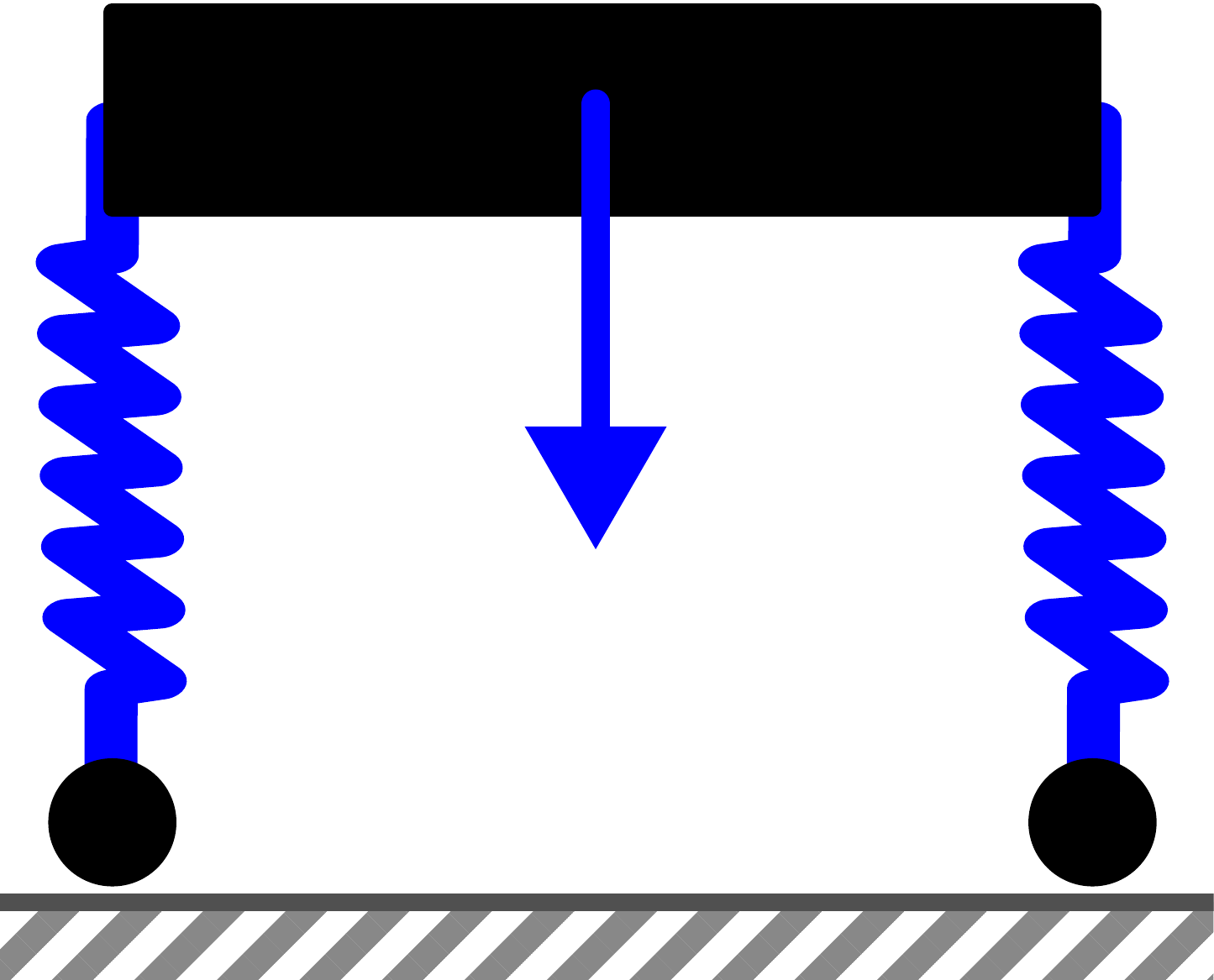}\\
\includegraphics[width=.9\columnwidth]{code/salt.pdf}
}
\caption{
\label{fig:ex}
Trajectory outcomes in mechanical systems subject to unilateral constraints.
(\emph{left})
In general, trajectory outcomes depend discontinuously on initial conditions.
In the pictured model for rigid--leg trotting (adapted from~\cite{RemyBuffinton2010}), discontinuities arise when two legs touch down:
if the legs impact simultaneously (corresponding to rotation $\theta(0) = 0$), then the post--impact rotational velocity is zero;
if the left leg impacts before the right leg ($\theta(0) > 0$, blue) or vice--versa ($\theta(0) < 0$, red), 
then the post--impact rotational velocities are bounded away from zero.
(\emph{right})
In the pictured model for soft--leg trotting 
(adapted from~\cite{BurdenGonzalezVasudevan2015tac} with the addition of a nonlinear damper coupling the body and limbs), 
trajectory outcomes (solid lines) are continuous and piecewise--differentiable at $\theta(0)=0$ (dashed lines). 
}
}
\end{figure}

To move through and interact with the world, terrestrial agents intermittently contact terrain and objects.
The dynamics of this interaction are, to a first approximation, \emph{hybrid}, with transitions between contact modes summarized by abrupt changes in system velocities~\cite{JohnsonBurden2016ijrr}.
Such phenomenological models are known in general to exhibit a range of pathologies that plague hybrid systems, including non--existence or non--uniqueness of trajectories~\cite{HurmuzluMarghitu1994, Stewart2000}~\cite[Sec.~5]{Ballard2000}, or discontinuous dependence of trajectory outcomes on initial conditions (i.e. states and parameters)~\cite{RemyBuffinton2010}~\cite[Sec.~7]{Ballard2000}; see~\fig{ex} (\emph{left}).
Although instances of these pathologies can occur in physical systems~\cite{HinrichsOestreich1997},
these occurrences are rare in everyday experience involving locomotion and manipulation with limbs.
Our view is that these pathologies lie chiefly in the modeling formalism, and can be effectively removed by appropriately restricting the models under consideration without loss of relevance for many physical systems of interest.

Specifically, this paper provides mathematical conditions on \emph{mechanical} systems subject to \emph{unilateral} constraints that ensure trajectory outcomes vary continuously and piecewise--differentiably with respect to initial conditions.
Conditions that ensure continuity are known; see for instance Schatzman's work on the one--dimensional impact problem~\cite{Schatzman1998} or Ballard's seminal result~\cite[Thm.~20]{Ballard2000}.
Furthermore, when the sequence of constraint activations and deactivations is held fixed, it has been known for some time that outcomes depend differentiably on initial conditions; see~\cite{AizermanGantmacher1958} for the earliest instance of this result we found in the English literature and~\cite{HiskensPai2000, GrizzleAbba2002, WendelAmes2012, BurdenRevzen2015tac} for modern treatments.
Our contribution is a proof that imposing an additional \emph{admissibility} condition ensures continuous trajectory outcomes are \emph{piecewise--}differentiable with respect to initial conditions, even as the sequence of constraint activations and deactivations varies; see~\fig{ex} (\emph{right}).
The operative notion of piecewise--differentiability was originally developed by the nonsmooth analysis community to study structural stability of nonlinear programs~\cite{Robinson1987}, and now forms the basis of a rather complete generalization of Calculus accommodating non--linear first--order approximations~\cite{Scholtes2012}.
In the terminology of that community, we provide conditions that ensure the flow of a mechanical system subject to unilateral constraints is $PC^r$, and therefore possesses a piecewise--linear \emph{Bouligand} (or \emph{B--})derivative.

As discussed in more detail in~\sct{disc}, we envision the existence and straightforward computability of the B--derivative of the flow to be useful in practice because it supports generalization of familiar control techniques to a class of hybrid systems with physical significance.
In particular, building on related work that dealt with differential equations with discontinuous right--hand--sides~\cite{BurdenSastry2016siads, Burden2014phd}, our B--derivative can be used to assess stability, controllability, or optimality of trajectories in mechanical systems subject to unilateral constraints.
As control of dynamic and dexterous robots increasingly relies on scalable algorithms for optimization and learning that presume the existence of first--order approximations (i.e. gradients or gradient--like objects)~\cite{MombaurLongman2005, KuindersmaDeits2015, LevineFinn2016, KumarTodorov2016}, it is important to place application of such algorithms on a firm theoretical foundation.
From a theoretical perspective, the results in this paper dovetail with recent advances in simulation of hybrid systems~\cite{BurdenGonzalezVasudevan2015tac} in that one of the conditions necessary for the B--derivative to exist (namely, continuity of trajectory outcomes) is also requisite for convergence of numerical simulations.
Taken together, these observations suggest that a unified analytical and computational framework for modeling and control of mechanical systems subject to unilateral constraints may be within reach.

\subsection{Organization}
\label{sec:org}
We begin in~\sct{mdl} by specifying the class of dynamical systems under consideration, namely, \emph{mechanical} systems subject to \emph{unilateral} constraints.
\sct{diffi} summarizes the well--known fact that, when the contact mode sequence is fixed, trajectories vary differentiably with respect to initial conditions.
In~\sct{(dis)cont}, we observe (as others have) that trajectories generally vary discontinuously with respect to initial conditions as the contact mode sequence varies, but provide a 
sufficient condition that is known to restore continuity.
\sct{diffa} leverages continuity to provide conditions under which trajectories vary piecewise--differentiably with respect to initial conditions across 
contact
mode sequences, and \sct{disc} discusses extensions and implications for a systems theory for mechanical systems subject to unilateral constraints.

\subsection{Relation to prior work}
\label{sec:prior}
The technical content in~\sct{mdl},~\sct{diffi}, and~\sct{(dis)cont} appeared previously in the literature and is (more--or--less) well--known; we collate the results here in a sequence of technical Lemmas%
\footnote{For uniformity and clarity of exposition, we present previous results here as Lemmas regardless of the form in which they originally appeared. 
}
to contextualize our contributions in~\sct{diffa}.

\section{Mechanical systems subject to unilateral constraints}
\label{sec:mdl}

In this paper, we study the dynamics of a mechanical system with configuration coordinates $q\in Q=\R^d$ subject to (perfect, holonomic, scleronomic)%
\footnote{A constraint is:
\emph{perfect} if it only generates force in the direction normal to the constraint surface; 
\emph{holonomic} if it varies with configuration but not velocity;
\emph{scleronomic} if it does not vary with time.
}
unilateral constraints $a(q) \ge 0$ 
specified by a differentiable function %
$a \colon Q\into \R^n$
where $d,n\in\N$ are finite.
We are primarily interested in systems with $n > 1$ constraints,  
whence we regard the inequality $a(q) \ge 0$ as being enforced componentwise.
Given any $J\subset\set{1,\dots,n}$, 
and letting $\abs{J}$ denote the number of elements in the set $J$,
we let 
$a_J \colon Q \into \R^{\abs{J}}$ 
denote the function obtained by selecting the component functions of $a$ indexed by $J$, 
and we regard the equality $a_J(q) = 0$ as being enforced componentwise.
It is well--known~\see{e.g.~\cite[Sec.~3]{Ballard2000} or \cite[Sec.~2.4,~2.5]{JohnsonBurden2016ijrr}}
that with
$J = \set{j\in\set{1,\dots,n} : a_j(q) = 0}$ 
the system's dynamics take the form
\begin{subequations}\label{eq:dyn}
\begin{align}
  M(q)\ddot{q} & = f(q,\dot{q}) + c(q,\dot{q})\dot{q} + Da_J(q)^\tr \lambda_J(q,\dot{q}),\label{eq:dyn:cont}\\
  \dot{q}^+ & = \Delta_J(q,\dot{q}^-)\dot{q}^-,\label{eq:dyn:disc}
\end{align}
\end{subequations}
where
$M \colon Q \into \R^{d\times d}$
specifies the mass matrix (or \emph{inertia tensor}) for the mechanical system in the $q$ coordinates,
$f \colon TQ \into \R^d$ 
is termed the \emph{effort map}~\cite{Ballard2000}
and specifies%
\footnote{We let $TQ = \R^d\times\R^d$ denote the \emph{tangent bundle} of the configuration space $Q$; an element $(q,\dot{q})\in TQ$ can be regarded as a pair containing a vector of generalized configurations $q\in\R^d$ and velocities $\dot{q}\in\R^d$; we write $\dot{q}\in T_q Q$.}
the internal and applied
forces, 
$c \colon TQ \into \R^{d\times d}$ 
denotes the \emph{Coriolis matrix} 
determined%
\footnote{For each $\ell,m\in\set{1,\dots,d}$
the $(\ell,m)$ entry $c_{\ell m}$ is determined 
from the entries of $M$ via
$c_{\ell m} = -\frac{1}{2} \sum_{k=1}^d \paren{ D_k M_{\ell m} + D_m M_{\ell k} - D_\ell M_{k m}}$.
}
by $M$,
$Da_J \colon Q \into \R^{\abs{J}\times d}$ 
denotes the (Jacobian) derivative of the constraint function $a_J$ with respect to the coordinates,
$\lambda_J \colon TQ \into \R^{\abs{J}}$ 
denotes the reaction forces generated in contact mode $J$ to enforce the constraint $a_J(q) \ge 0$,
$\Delta_J \colon TQ \into \R^{d\times d}$ 
specifies the collision restitution law that instantaneously resets velocities to ensure compatibility with the constraint $a_J(q) = 0$,
\eqnn{
\Delta_J(q,\dot{q}) =& \I - (1+\gamma(q,\dot{q})) M(q)^{-1} Da_J(q)^\tr \Lambda_J(q) Da_J(q),
}
where 
$\I$ is the $d$--dimensional identity matrix,
$\gamma \colon TQ\into[0,\infty)$ specifies the \emph{coefficient of restitution},
$\dot{q}^+$ (resp. $\dot{q}^-$) denotes the right-- (resp. left--)handed limits of the velocity vector with respect to time,
and $\Lambda_J\colon Q\into R^{d\times d}$
\eqnn{
\Lambda_J(q) = \paren{Da_J(q) M(q)^{-1} Da_J(q)^\tr}^{-1}.
}

\defna{contact modes}{
\label{def:contact}
With \\$A =\set{q\in Q : a(q) \ge 0}$ denoting the set of \emph{admissible} configurations,
the constraint functions $\set{a_j}_{j=1}^n$ partition $A$ into a finite collection%
\footnote{We let $2^n = \set{J \subset \set{1,\dots,n}}$ denote the \emph{power set} (i.e. the set containing all subsets) of $\set{1,\dots,n}$.}
$\set{A_J}_{J\in 2^n}$ 
of \emph{contact modes}:
\eqnn{
\forall J\in 2^n : A_J = \left\{q\in Q \mid \right. & a_J(q) = 0,\\ 
&\left. \forall i\not\in J : a_i(q) > 0\right\}.
}
We let $TA = \set{(q,\dot{q})\in TQ : q\in A}$ and \\$TA_J = \set{(q,\dot{q})\in TQ : q\in A_J}$ for each $J\in 2^n$.
}

\rem{
\label{rem:contact}
In~{\defcontact}, 
$J = \set{1,\dots,n}$ indexes the maximally constrained contact mode and 
$J = \emptyset$ indexes the unconstrained contact mode. 
Since any velocity is allowable in the unconstrained mode, we adopt the convention 
$\Delta_\emptyset(q,\dot{q})=\I$.
}
In the present paper, we will assume that appropriate conditions have been imposed to ensure trajectories of~\eqref{eq:dyn} exist on a region of interest in time and state.

\asmpna{existence and uniqueness}{
\label{asmp:flow}
There exists a \emph{flow} for~\eqref{eq:dyn},
that is,
a function $\phi:\e{F}\into TA$
where 
$\e{F}\subset[0,\infty)\times TA$ is an open subset 
containing $\set{0}\times TA$ 
and for each $(t,(q,\dot{q}))\in\e{F}$
the restriction
$\phi|_{[0,t]\times\set{(q,\dot{q})}}:[0,t]\into TQ$
is the unique left--continuous trajectory for~\eqref{eq:dyn} initialized at $(q,\dot{q})$.
}

\rem{
The problem of ensuring trajectories of~\eqref{eq:dyn} exist and are unique has been studied extensively; we refer the reader to~\cite[Thm.~10]{Ballard2000} for a specific result and~\cite{JohnsonBurden2016ijrr} for a general discussion of this problem.
}

Since we are concerned with differentiability properties of the flow, 
we assume the elements in~\eqref{eq:dyn} are differentiable.

\asmpna{$C^r$ vector field and reset map}{
\label{asmp:diff}
The vector field~\eqref{eq:dyn:cont} and reset map~\eqref{eq:dyn:disc} are continuously differentiable to order $r\in\N$.
}

\rem{
If we restricted our attention to the continuous--time dynamics in~\eqref{eq:dyn}, then Assump.~\ref{asmp:diff} would suffice to provide the local existence and uniqueness of trajectories imposed by Assump.~\ref{asmp:flow};
as illustrated by~\cite[Ex.~2]{Ballard2000},
Assump.~\ref{asmp:diff} does not suffice when the vector field~\eqref{eq:dyn:cont} is coupled to the reset map~\eqref{eq:dyn:disc}.
}

\section{Differentiability within contact mode sequences}
\label{sec:diffi}

It is possible to satisfy {\asmpflow} under mild conditions that allow trajectories to 
exhibit phenomena such as \emph{grazing} (wherein the trajectory activates a new constraint without undergoing impact) or \emph{Zeno} (wherein the trajectory undergoes an infinite number of impacts in a finite time interval).
In this and subsequent sections, where we seek to study differentiability properties of the flow, we will not be able to accommodate grazing or Zeno phenomena.
Therefore we proceed to restrict the trajectories under consideration.

\defna{constraint activation/deactivation}{
\label{def:act}
The trajectory $\phi^{(q,\dot{q})}$ initialized at $(q,\dot{q})\in TA_J\subset TQ$
\emph{activates constraints $I\in 2^n$ at time $t > 0$} if
(i) no constraint in $I$ was active immediately before time $t$
and
(ii) all constraints in $I$ become active at time $t$.
Formally,%
\footnote{$\phi( (t_1,t_2), (q,\dot{q}) ) = \set{\phi(t,(q,\dot{q})) : t\in (t_1,t_2)}\subset TQ$ denotes the \emph{image} of $\phi^{(q,\dot{q})}$ over the interval $(t_1,t_2)\subset[0,\infty)$.}
\eqnn{
\exists\veps > 0 :\ I\cap J = \emptyset,\ \text{(i)}\ &\phi\paren{(t-\veps,t),(q,\dot{q})}\subset TA_J,\\
\text{(ii)}\ &\phi(t,(q,\dot{q}))\in TA_{I\cup J}.
}
We refer to $t$ as a \emph{constraint activation time} for $\phi^{(q,\dot{q})}$.
Similarly, the trajectory $\phi^{(q,\dot{q})}$ 
\emph{deactivates constraints $I\in 2^n$ at time $t > 0$} if
(i) all constraints in $I$ were active at time $t$
and
(ii) no constraint in $I$ remains active immediately after time $t$.
Formally,
\eqnn{
\exists\veps > 0 :\ I\subset J,\ \text{(i)}\ &\phi(t,(q,\dot{q}))\in TA_{J},\\
\text{(ii)}\ &\phi\paren{(t,t+\veps),(q,\dot{q})}\subset TA_{J\sm I}.
}
We refer to $t$ as a \emph{constraint deactivation time} for $\phi^{(q,\dot{q})}$.
}

\defna{admissible activation/deactivation}{
\label{def:adact}
A constraint activation time $t > 0$ for $\phi^{(q,\dot{q})}$ is \emph{admissible} if the constraint velocity%
\footnote{Formally, the \emph{Lie derivative}~\cite[Prop.~12.32]{Lee2012} of the constraint along the vector field specified by~\eqref{eq:dyn:cont}.
Although constraint functions are technically only functions of configuration $q\in Q$ and not the full state $(q,\dot{q})\in TQ$, by a mild abuse of notation we allow ourselves to consider compositions $a\circ\phi$ rather than the formally correct $a\circ\pi_Q\circ\phi$ where $\pi_Q\colon TQ\into Q$ is the canonical projection.}
for all activated constraints $I\in 2^n$ is negative.
Formally, with $(\rho,\dot{\rho}^-) = \lim_{s\goesto t^-}\phi(s,(q,\dot{q}))$ denoting the left--handed limit of the trajectory at time~$t$,
\eqnn{
\forall i\in I: D_t\brak{a_i\circ\phi}(0,(\rho,\dot{\rho}^-)) = Da_i(\rho)\dot{\rho}^- < 0.
}
A constraint deactivation time $t > 0$ for $\phi^{(q,\dot{q})}$ is \emph{admissible} if, for all deactivated constraints $I\in 2^n$: 
(i) the constraint velocity or constraint acceleration%
\footnote{Formally, the second Lie derivative of the constraint along the vector field specified by~\eqref{eq:dyn:cont}.}
is positive,
or
(ii) the time derivative of the contact force is negative.
Formally, with $(\rho,\dot{\rho}^+) = \lim_{s\goesto t^+}\phi(s,(q,\dot{q}))$ denoting the right--handed limit of the trajectory at time $t$,
for all $i\in I:$
\eqnn{\label{eq:deact}
\text{(i)}\ D_t\brak{a_i\circ\phi}(0,(\rho,\dot{\rho}^+)) &> 0\ \text{or}\\
D^2_t\brak{a_i\circ\phi}(0,(\rho,\dot{\rho}^+)) &> 0,\\ 
\text{or (ii)}\ D_t\brak{\lambda_i\circ\phi}(0,(\rho,\dot{\rho}^+)) &< 0. 
}
}

\rem{\label{rem:adact}
The conditions for admissible constraint deactivation in case (i) of~\eqref{eq:deact} can only arise at admissible constraint activation times; otherwise the trajectory is continuous, whence active constraint velocities and accelerations are zero.
}

\defna{admissible trajectory}{
\label{def:trj}
A trajectory $\phi^{(q,\dot{q})}$ is \emph{admissible on $[0,t]\subset\R$} if 
(i) it has a finite number of constraint activation (hence, deactivation) times on $[0,t]$, 
and 
(ii) every constraint activation and deactivation is admissible;
otherwise the trajectory is \emph{inadmissible}.
}

\remna{admissible trajectories}{
\label{rem:trj}
The key property admissible trajectories possess that will be leveraged in what follows is: 
time--to--activation and time--to--deactivation are differentiable with respect to initial conditions;
the same is not generally true of inadmissible trajectories.
}

\remna{grazing is not admissible}{
\label{rem:graze}
The restriction in~{\deftrj} 
that all constraint activation/deactivation times are admissible precludes admissibility of grazing.
}

\remna{Zeno is not admissible}{
\label{rem:zeno}
The restriction in~{\deftrj} 
that a finite number of constraint activations 
occur on a compact time interval precludes admissibility of Zeno.
}

\defna{contact mode sequence}{
\label{def:seq}
The \emph{contact mode sequence}%
\footnote{This definition differs from the \emph{word} of~\cite[Def.~4]{JohnsonBurden2016ijrr} in that a contact mode is included in the sequence only if nonzero time is spent in the mode; this definition is more closely related to the \emph{word} of~\cite[Eqn.~72]{BurdenSastry2016siads}} 
associated with a trajectory $\phi^{(q,\dot{q})}$ that is admissible on $[0,t]\subset\R$ 
is the unique function $\omega\colon\set{0,\dots,m}\into 2^n$ such that there exists a finite sequence of times $\set{t_\ell}_{\ell=0}^{m+1}\subset [0,t]$ for which $0 = t_0 < t_1 < \cdots < t_{m+1} = t$ and
\eqnn{
\forall \ell\in\set{0,\dots,m} : \phi((t_{\ell},t_{\ell+1}),(q,\dot{q})) \subset TA_{\omega(\ell)}.
}
}

\rem{
\label{def:seq}
In~{\defseq}, the sequence $\omega$ is easily seen to be unique by the admissibility of the trajectory; indeed, the associated time sequence consists of start, stop, and constraint activation/deactivation times.
}

\asmpna{independent constraints}{
\label{asmp:indep}
The constraints are independent:
\eqnn{
\forall J\in 2^n,\ q\in a_J^{-1}(0) : &\set{Da_j(q)}_{j\in J}\subset T^*_q Q\\
&\quad\text{is linearly independent}.
}
}

\rem{
\label{rem:indep}
Algebraically, {\asmpindep} 
implies 
that the constraint forces $\lambda_J$ are well--defined,
and
that there are no more constraints than degrees--of--freedom, $n \le d$.
Geometrically, it implies for each $J\in 2^n$ that $a_J^{-1}(0)\subset Q$ is an embedded codimension--$\abs{J}$ submanifold, and that the codimension--1 submanifolds $\set{a_j^{-1}(0)}_{j\in J}$ intersect transversally; this follows from~\cite[Thm.~5.12]{Lee2012} since each $a_J\colon Q\into\R$ must be constant--rank on its zero section.
}

We now state the well--known fact%
\footnote{%
The result follows via a straightforward composition of smooth flows with smooth time--to--impact maps; we refer the interested reader to~\cite[App.~A1]{BurdenRevzen2015tac} for details.
}
that, if the contact mode sequence is fixed, then admissible trajectory outcomes are differentiable with respect to initial conditions.

\lemna{differentiability within mode seq.~\cite{AizermanGantmacher1958}}{
\label{lem:diffi}
Under~{\asmpflow},~{\asmpdiff}, and~{\asmpindep}, 
with $\phi\colon[0,\infty)\times TA\into TA$ denoting the flow,
if $\Sigma\subset TQ$ is a $C^r$ embedded submanifold such that all trajectories initialized in $\Sigma\cap TA$ 
\begin{enumerate}
  \item[(i)] are admissible on $[0,t]\subset\R$
  and
  \item[(ii)] have the same contact mode sequence,
\end{enumerate}
then 
the restriction
$\phi|_{\set{t}\times\Sigma}$
is $C^r$.
}

\section{(Dis)continuity across contact mode sequences}
\label{sec:(dis)cont}

As stated in~\sct{intro}, the point of this paper is to 
provide sufficient conditions that ensure trajectories of~\eqref{eq:dyn} vary differentiably as the contact mode sequence varies.
A precondition for differentiability is continuity, whence in this section we consider what condition must be imposed to give rise to continuity in general.
We begin in~\sct{discont} by demonstrating that the transversality of constraints imposed by~{\asmpindep} generally gives rise to discontinuity, then introduce an orthogonality condition in~\sct{cont} that suffices to restore continuity.

\subsection{Discontinuity across contact mode sequences}
\label{sec:discont}

Consider an unconstrained initial condition $(q,\dot{q})\in TA_\emptyset\subset TQ$ that impacts (i.e. admissibly activates) exactly two constraints $i,j\in\set{1,\dots,n}$ at time $t > 0$;  
with $(\rho,\dot{\rho}^-) = \phi(t,(q,\dot{q}))$ we have
\eqnn{
a_{\set{i,j}}(\rho) = 0,\
Da_i(\rho)\dot{\rho}^- < 0,\ 
Da_j(\rho)\dot{\rho}^- < 0.\\
}
The pre--impact velocity $\dot{\rho}^-$ abruptly resets via~\eqref{eq:dyn:disc}:
\eqnn{
\label{eq:ij}
\dot{\rho}^+ = \Delta_{\set{i,j}}(\rho)\dot{\rho}^-.
}
As noted in~{\remindep},
the constraint surfaces $a_i^{-1}(0)$, $a_j^{-1}(0)$ intersect transversally.
Therefore given any $\veps > 0$ it is possible to find $(q_i,\dot{q}_i)$ and $(q_j,\dot{q}_j)$ in the open ball of radius $\veps$ centered at $(q,\dot{q})$
such that 
the trajectory $\phi^{(q_i,\dot{q}_i)}$ impacts constraint $i$ before constraint $j$ 
and
$\phi^{(q_j,\dot{q}_j)}$ impacts $j$ before $i$.
As $\veps > 0$ tends toward zero, the time spent flowing according to~\eqref{eq:dyn:cont} tends toward zero, hence the post--impact velocities tend toward the twofold iteration of~\eqref{eq:dyn:disc}:
\eqnn{
\label{eq:ijji}
\dot{\rho}_i^+ = \Delta_{\set{i,j}}(\rho)\Delta_i(\rho)\dot{\rho}^-,\\ 
\dot{\rho}_j^+ = \Delta_{\set{i,j}}(\rho)\Delta_j(\rho)\dot{\rho}^-.
}
Recalling for all $J\in 2^n$ that $\Delta_J\in\R^{d\times d}$ is an orthogonal projection%
\footnote{relative to the inner product $\ipM{\cdot}{\cdot}$}
onto the tangent plane of the codimension--$\abs{J}$ submanifold $a_J^{-1}(0)$,
observe that $\dot{\rho}_i^+ = \dot{\rho}^+ = \dot{\rho}_j^+$ if and only if $Da_i(\rho)$ is orthogonal to $Da_j(\rho)$.
Therefore if constraints intersect transversally but non--orthogonally, outcomes from the dynamics in~\eqref{eq:dyn} vary discontinuously as the contact mode sequence varies.

\remna{discontinuous locomotion outcomes}{
The analysis of a saggital--plane quadruped in~\cite{RemyBuffinton2010} provides an instructive example of the behavioral consequences of transverse but non--orthogonal constraints in a model of legged locomotion.
As summarized in~\cite[Table~2]{RemyBuffinton2010}, the model possesses 3 distinct but nearby trot (or trot--like) gaits, corresponding to whether two legs impact simultaneously (as in~\eqref{eq:ij}) or at different time instants (as in~\eqref{eq:ijji});
the trot that undergoes simultaneous impact is unstable due to discontinuous dependence of trajectory outcomes on initial conditions.
}

\subsection{Continuity across contact mode sequences}
\label{sec:cont}

To preclude the discontinuous dependence on initial conditions exhibited in~\sct{discont}, we strengthen the \emph{transversality} of constraints implied by~{\asmpindep} by imposing \emph{orthogonality} of constraints.

\asmpna{orthogonal constraints}{
\label{asmp:orth}
Constraint surfaces intersect orthogonally:
\eqnn{\label{eq:orth}
\forall i,j\in\set{1,\dots,n},\ i\ne j,\ q\in a_i^{-1}(0)\cap a_j^{-1}(0) :\quad\\
\ip{Da_i(q)}{Da_j(q)}_{M^{-1}} = 0.
}
}

\rem{
\label{rem:orth}
Note that~{\asmporth} is strictly stronger than~{\asmpindep}.
Physically, the assumption can be interpreted as asserting that any two independent limbs that can undergo impact simultaneously must be inertially decoupled.
This can be achieved in artifacts by introducing series compliance in a sufficient number of degrees--of--freedom.
}

\sct{discont} demonstrated that~{\asmporth} is necessary in general to preclude discontinuous dependence on initial conditions.
The following result demonstrates that this assumption is sufficient to ensure continuous dependence on initial conditions.

\lemna{continuity across mode seq.~{\cite[Thm.~20]{Ballard2000}}}{
\label{lem:cont}
Under~{\asmpflow},~{\asmpdiff}, and {\asmporth},
with $\phi\colon[0,\infty)\times TA\into TA$ denoting the flow,
if 
$t\in\R$
and
$(\q,\dot{\q})\in TA\subset TQ$ 
are such that
$t$ is not a constraint activation time for $(\q,\dot{\q})$,
then $\phi$ is continuous at $(t,(\q,\dot{\q}))$.
}

\remna{continuity across mode seq.}{
The preceding result implies that the flow $\phi$ is continuous almost everywhere in both time and state, without needing to restrict to admissible trajectories.
Thus orthogonal constraints ensure the flow $\phi$ depends continuously on initial conditions, even along trajectories that exhibit grazing and Zeno phenomena.%
\footnote{%
We remark that~\cite[Thm.~20]{Ballard2000} implies the function $\phi$ is continuous everywhere with respect to the quotient metric defined in~\cite[Sec.~III]{BurdenGonzalezVasudevan2015tac}, whence the numerical simulation algorithm in~\cite[Sec.~IV]{BurdenGonzalezVasudevan2015tac} is provably--convergent for all trajectories (even those that graze) up to the first occurrence of Zeno.
}
For the reason described in~{\remtrj},
we will not be able to accommodate these phenomena when we study differentiability properties of trajectories in the next section.
}

\section{Differentiability across contact mode sequences}
\label{sec:diffa}

We now provide conditions that ensure
trajectories depend \emph{differentiably} on initial conditions, even as the contact mode sequence varies.
In general, 
the flow does not possess a classical Jacobian (alternately called \emph{\Frechet} or \mbox{\emph{F--})derivative}, i.e. there does not exist a single linear map that provides a first--order approximation for the flow.
Instead,
under the admissibility conditions introduced in~\sct{diffi},
we show that the flow admits a \emph{piecewise--}linear first--order approximation termed%
\footnote{This terminology was introduced, to the best of our knowledge, by Robinson~\cite{Robinson1987}.}
a \emph{Bouligand} (or \emph{B--})derivative~\cite[Ch.~3.1]{Scholtes2012}.
Though perhaps unfamiliar, this derivative is nevertheless quite useful.
Significantly, unlike functions that are merely directionally differentiable, B--differentiable functions admit generalizations of many techniques familiar from calculus, including
the Chain Rule~\cite[Thm~3.1.1]{Scholtes2012} (and hence Product and Quotient Rules~\cite[Cor.~3.1.1]{Scholtes2012}),
Fundamental Theorem of Calculus~\cite[Prop.~3.1.1]{Scholtes2012},
and Implicit Function Theorem~\cite[Thm.~4.2.3]{Scholtes2012},
and the B--derivative can be employed to implement scalable algorithms~\cite{Kiwiel1985} for optimization or learning.

We proceed by showing that the flow is piecewise--differentiable in the sense defined in~\cite[Ch.~4.1]{Scholtes2012} and recapitulated here;
functions that are piecewise--differentiable in this sense are always B--differentiable~\cite[Prop.~4.1.3]{Scholtes2012}.
Let $r\in\N\cup\set{\infty}$ denote an order of differentiability%
\footnote{We let context specify whether $r = \infty$ indicates ``mere'' smoothness or the more stringent condition of analyticity.}
and $D\subset\R^m$ be open.
A continuous function $\psi\colon D\into\R^\ell$ is called \emph{piecewise--$C^r$} if
the graph of $\psi$ is everywhere locally covered by the graphs of a finite collection of 
functions that are $r$ times continuously differentiable ($C^r$--functions).%
\footnote{The definition of piecewise--$C^r$ may at first appear unrelated to the intuition that a function ought to be piecewise--differentiable precisely if its ``domain can be partitioned locally into a finite number of regions relative to which smoothness holds''~\cite[Section~1]{Rockafellar2003}.
However, as shown in~\cite[Theorem~2]{Rockafellar2003}, piecewise--$C^r$ functions are always piecewise--differentiable in this intuitive sense.
}
Formally,
for every $x\in D$ there must exist an open set $U\subset D$ containing $x$ and a finite collection $\set{\psi_\omega\colon U\into\R^\ell}_{\omega\in\Omega}$ of $C^r$--functions such that for all $x\in U$
we have $\psi(x)\in\set{\psi_\omega(x)}_{\omega\in\Omega}$.

\begin{figure}[th]
\centering
{
\iftoggle{acm}
{
\def\svgwidth{1.5\columnwidth} 
\resizebox{\columnwidth}{!}{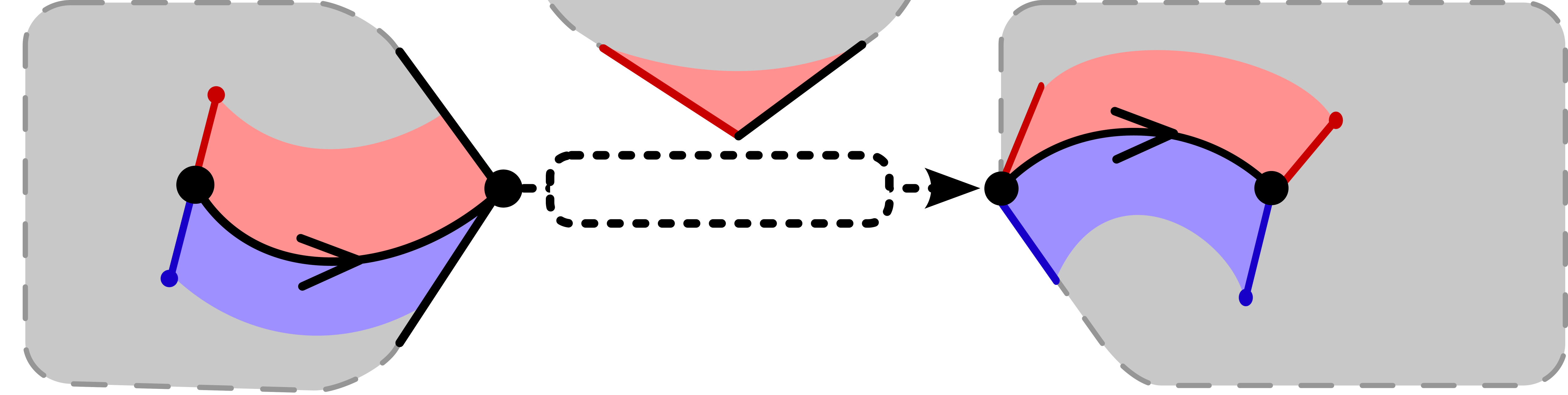}
}
{
\def\svgwidth{0.95\columnwidth} 
\resizebox{\columnwidth}{!}{\input{ad.pdf_tex}}
}
}
\caption{\label{fig:hds}
Illustration of trajectory undergoing simultaneous
constraint activation and deactivation:
the trajectory initialized at $(\q,\dot{\q})\in TA_{\set{1}}\subset TQ$
flows via~\eqref{eq:dyn:cont} to a point $(\rho,\dot{\rho}^-)\in TA_{\set{1}}$
where both 
the constraint force $\lambda_1$ and constraint function $a_2$ are zero,
instantaneously resets velocity via~\eqref{eq:dyn:disc} to $\dot{\rho}^+ = \Delta_{\set{2}}(\rho)\dot{\rho}^-$,
then flows via~\eqref{eq:dyn:cont} to $\phi(t,(\q,\dot{\q}))\in TA_{\set{2}}\subset TQ$.
Nearby trajectories undergo activation and deactivation at distinct times: 
trajectories initialized in the red region, e.g. $(v_r,\dot{v}_r)$, deactivate constraint $1$ and flow through contact mode $TA_\emptyset$ before activating constraint $2$---their contact mode sequence is $\paren{\set{1},\emptyset,\set{2}}$---while
trajectories initialized in the blue region, e.g. $(v_b, \dot{v}_b)$, activate $2$ and flow through $TA_{\set{1,2}}$ before deactivating $1$---their contact mode sequence is $\paren{\set{1},\set{1,2},\set{2}}$.
Piecewise--differentiability of the trajectory outcome is illustrated by the fact that red outcomes lie along a different subspace than blue.
}
\end{figure}

We now state and prove the main result of this paper:  whenever the flow of a mechanical system subject to unilateral constraints is continuous and admissible, it is piecewise--$C^r$; see~\fig{hds} for an illustration.

\thmna{piecewise--differentiable flow}{
\label{thm:diffb}
Under~{\asmpflow},~{\asmpdiff}, and~{\asmporth},
with $\phi\colon [0,\infty)\times TA\into TA$ denoting the flow,
if 
$t\in [0,\infty)$,
$(\q,\dot{\q})\in TA\subset TQ$,
and
$\Sigma\subset TQ$ is a $C^r$ embedded submanifold containing $(\q,\dot{\q})$
such that 
\begin{enumerate}
\item[(i)] the trajectory $\phi^{(\q,\dot{\q})}$ activates and/or deactivates constraints at time $s\in(0,t)$,
\item[(ii)] $\phi^{(\q,\dot{\q})}$ has no other activation or deactivation times in $[0,t]$,
\item[(iii)] trajectories initialized in $\Sigma\cap TA$ are admissible on $[0,t]$, 
and 
\item[(iv)] the set $\Omega$ of contact mode sequences for trajectories initialized in $\Sigma\cap TA$ is finite,
\end{enumerate}
then 
the restriction
$\phi|_{[0,\infty)\times\Sigma}$
is piecewise--$C^r$ at $(t,(\q,\dot{\q}))$.
}
\begin{proof}
We seek to show that 
the restriction
$\phi|_{[0,\infty)\times\Sigma}$
is piecewise--$C^r$ at $(t,(\q,\dot{\q}))$.
We will proceed by constructing a finite set of $r$ times continuously differentiable selection functions for $\phi$ on $[0,t]\times \Sigma$.
In the example given in \figref{fig:hds}, there are two selection functions, one corresponding to a 
perturbation along $(v_r,\dot{v}_r)$, colored red, and the other along $(v_b,\dot{v}_b)$, colored blue.
These selection functions will be indexed by a pair of functions $(\omega,\eta)$ where: 
$\omega\colon \set{0,\dots,m}\into 2^n$ is a contact mode sequence, 
i.e. $\omega\in\Omega$;
$\eta\colon \set{0,\dots,m-1}\into\set{1,\dots,n}$ indexes constraints that undergoes admissible activation or deactivation%
\footnote{%
In light of {\remadact}, we only consider deactivations of type (ii) in {\defadact}.
In some systems, a deactivation of type (ii) may only arise following a (simultaneous) activation; it suffices to restrict to functions $\eta$ that do not index such deactivations.
}
at the contact mode transition indexed by $\ell\in\set{0,\dots,m-1}$. 
For instance, in \figref{fig:hds} the index functions for the (de)activation sequence starting from
$(v_r,\dot{v}_r)$, in red,
are $\w_r(0)= \set{1}, \w_r(1)=\emptyset, \w_r(2)=\set{2}$, $\eta_r(1)=1, \eta_r(2)=2$,
and the index functions for the (de)activation sequence starting from $(v_b, \dot{v}_b)$, in blue,
are $\w_b(0)= \set{1}, \w_b(1)=\set{1,2}, \w_b(2)=\set{2}$, $\eta_b(1)=2, \eta_b(2)=1$.
Note that for each $\omega\in\Omega$ the set $H(\omega)$ of possible $\eta$'s is finite; since the set $\Omega$ is finite by assumption, the set of pairs $(\omega,\eta)$ is finite.

Let $(\omega\colon \set{0,\dots,m}\into 2^n)\in\Omega$ and $(\eta\colon \set{0,\dots,m-1}\into\set{1,\dots,n})\in H(\omega)$ be as described above. 
Let $\mu \colon \set{0,\dots,m}\into 2^n$ be defined as 
$\mu(k) = \bigcup_{i=0}^{k-1}\set{\eta(i)}$, where we adopt the convention that $\bigcup_{i=0}^{-1}\set{i} = \emptyset$;
note that $\mu$ is uniquely determined by $\eta$.%
\footnote{$\eta$ is not uniquely determined by $\w$ due to the possibility of instantaneous activation/deactivation for the same constraint; consider for instance the bounce of an elastic ball~\cite[Ch.~2.4]{GuckenheimerHolmes1983}.}
For the sake of readability, we suppress dependence on $\eta$ and $\w$ until~\eqref{eq:phiomegaeta}. 
Let $(\rho,\dot{\rho}^-)=\lim_{u\uparrow s}\phi(u,(\q,\dot{\q}))$.
For all $k\in\set{0,\dots,m}$
define $\dot{\rho}_{k} = \Delta_{\mu(k)}(\rho)\dot{\rho}^-$.
There exists an open neighborhood $U_{k} \subset TQ$ containing $(\rho,\dot{\rho}_{k})$ such that the vector field determined by~\eqref{eq:dyn:cont} at $\w(k)$ admits a $C^r$ extension to $F_{k}\colon U_{k}\into\R^{2d}$.
(Note that for $k= m$ (resp. $k = 0$) the neighborhood $U_{k}$ can be taken to additionally include $\phi( (s,t], (\q,\dot{\q}) )$ (resp. $\phi( [0,s), (\q,\dot{\q}) )$).)

By the Fundamental Theorem on Flows~\cite[Thm.~9.12]{Lee2012}, $F_{k}$ determines a unique maximal flow $\phi_{k}\colon \e{F}_{k}\into U_{k}$ over a maximal flow domain $\e{F}_{k}\subset\R\times U_{k}$, which is an open set that contains $\set{0}\times U_{k}$, 
and the flow $\phi_{\ell}$ is $C^r$.
(Note that $(t-s,(\rho,\dot{\rho}_{m}))\in\e{F}_{m}$ and $(s,(\q,\dot{\q}))\in\e{F}_{0}$.)

If $\eta(\ell)$ indexes an admissible constraint activation (recall that $\ell\in\set{0, \dots, m-1}$),
then
there exists a time--to--activation 
$\tau_{\ell}\colon U_{\ell}\into\R$
defined over an open set
$U_{\ell}\subset TQ$ containing $(\rho,\dot{\rho}_{\ell})$ such that
\eqnn{
\forall (\p,\dot{\p})\in U_{\ell}: 
a_{\eta(\ell)}\circ\phi_{\ell}(\tau_{\ell}(\p,\dot{\p}),(\p,\dot{\p})) = 0.
}
If instead $\eta(\ell)$ indexes an admissible constraint deactivation,
then
there exists a time--to--deactivation $\tau_{\ell}\colon U_{\ell}\into\R$
defined over an open set
$U_{\ell}\subset TQ$ containing $(\rho,\dot{\rho}_{\ell})$ such that
\eqnn{
\forall (\p,\dot{\p})\in U_{\ell} : 
\lambda_{\eta(\ell)}\circ\phi_{\ell}(\tau_{\ell}(\p,\dot{\p}),(\p,\dot{\p})) = 0.
}
In either case, $\tau_{\ell}$ exists and is $C^r$ by the Implicit Function Theorem~\cite[Thm.~C.40]{Lee2012} due to admissibility of trajectories initialized in $\Sigma$.
(Note for $\ell = 0$ the neighborhood $U_{\ell}$ can be extended to include $\phi( [0,s), (\q,\dot{\q}) )$ using the semi--group property%
\footnote{$\phi_{\ell}(u + v,x) = \phi_{\ell}(u,\phi_{\ell}(v,x))$ whenever $(v,x), (u+v,x), (u,\phi_{\ell}(v,x))\in\e{F}_{\ell}$.}
of the flow $\phi_{\ell}$.)
See~\fig{hds} for an illustration of constraint activations and deactivations.

Let $\vphi_{\ell}\colon \R\times U_{\ell}\into \R\times U_{\ell}$ be defined for all $(u,(\p,\dot{\p}))\in\R\times U_{\ell}$ by
\eqnn{
\vphi_{\ell}(u,(\p,\dot{\p})) = \paren{u - \tau_{\ell}(\p,\dot{\p}),\phi_{\ell}\paren{\tau_{\ell}(\p,\dot{\p}),(\p,\dot{\p})}}.
}
The map $\vphi_{\ell}$ flows a state $(\p,\dot{\p})$ using the vector field from contact mode $\omega(\ell)$ until constraint $\eta(\ell)$ undergoes admissible activation/deactivation, and deducts the time required from the given budget $u$.
The total 
derivative of $\vphi_\ell$ at $(0, (\rho, \dot{\rho}_\ell))$ (see also \mbox{\cite[\S~7.1.4]{BurdenSastry2016siads}}) is
\eqnn{
D\vphi_{\ell}(0,(\rho,\dot{\rho}_\ell)) = \mat{cc}{
	1 & \frac{1}{gf}g \\
	0 & \II-\frac{1}{gf}fg
},
}
where $f= F(\rho, \dot{\rho}_\ell)$ and $g=Dh_{\eta(\ell)}(\rho, \dot{\rho}_\ell)$
where $h_\ell\colon TQ\into \R$ is defined for all $(q,\dot{q})\in TQ$ by $h_\ell(q,\dot{q})=a_{\eta(\ell)}(q)$. 

Let $\Gamma_{\ell}\colon \R\times TQ\into \R\times TQ$ be defined for all $(u,(\p,\dot{\p}))\in\R\times TQ$ by
\eqnn{
\Gamma_{\ell}(u,(\p,\dot{\p})) = (u,(\p,\Delta_{\mu(\ell)}(\p)\dot{\p})).
}
The map $\Gamma_{\ell}$ resets velocities to be compatible with contact mode $\omega(\ell)$ while leaving positions and times unaffected.
The total derivative of $\Gamma_\ell$ at $(u, (q, \dot{q}))$ is given by
\eqnn{
D\Gamma_{\ell}(u,(\p,\dot{\p}))  = \mat{ccc}{
	1 & 0 & 0 \\
	0 & \I & 0 \\
	0 & D_q(\Delta_{\mu(\ell)}(q)\dot{q}) & \Delta_{\mu(\ell)}(q) }.
}

For each $\omega\in\Omega$ and $\eta\in H(\omega)$ define $\phi_\omega^\eta$ by the formal composition
\eqnn{\label{eq:phiomegaeta}
\phi_\omega^\eta = \phi_{m}\circ\prod_{\ell=0}^{m-1}\paren{\Gamma_{\ell+1}\circ\vphi_{\ell}}.
}
We take as the domain of $\phi_\omega^\eta$ the set 
\eqnn{\e{F}_\omega^\eta = \paren{\phi_\omega^\eta}^{-1}( TQ )\subset\R\times TQ,}
noting that $\e{F}_\omega^\eta$ is 
(i) open since each function in the composition is continuous, and
(ii) nonempty since $(t,(\q,\dot{\q}))\in\e{F}_\omega^\eta$. 
The map $\phi_\omega^\eta$ flows states via a given contact mode sequence for a specified amount of time; note that some of the resulting ``trajectories'' are not physically realizable, as they may evaluate the flows $\set{\phi_k}_{k=0}^m$ in backward time.
An example of such a physically unrealizable ``trajectory'' is illustrated in~\figref{fig:hds} by 
$\phi_{\eta_r}^{\w_r}(t,(v_b,\dot{v}_b))$,
which first flows forward in time via the extended vector field $F_{\set{1}}$ past the constraint surface 
$\set{a_2(q)=0}$ until constraint 1 deactivates 
and then flows backwards in time until constraint 2 activates, ultimately terminating in $TA_{\set{2}}$.

With $\e{F} = \bigcap\set{\e{F}_\omega^\eta : \omega\in\Omega, \eta\in H(\omega)}\subset[0,\infty)\times TQ$,
for any $(u,(\p,\dot{\p}))\in \e{F}\cap\paren{[0,\infty)\times TA}$
with contact mode sequence $\omega\in\Omega$ 
and 
constraint sequence $\eta\in H(\omega)$, 
the trajectory outcome is obtained by applying $\phi_\omega^\eta$ to $(u,(\p,\dot{\p}))$,
i.e. $\phi(u,(\p,\dot{\p})) = \phi_\omega^\eta(u,(\p,\dot{\p}))$.
See~\fig{hds} for an illustration of trajectories with different contact mode sequences.

The maps 
$\vphi_{\ell}$,
$\Gamma_{\ell}$,
and
$\phi_\omega^\eta$
are $C^r$ on their domains since they are each obtained from a finite composition of $C^r$ functions.
Therefore the restriction%
\footnote{%
As a technical aside, we remark that the domain of $\phi$ is confined to $[0,\infty)\times TA$, 
whence invoking the definition of piecewise--differentiability 
requires a continuous extension $\td{\phi}$ of $\phi$ defined on a neighborhood of $(t,(\q,\dot{\q}))$ that is open relative to $[0,\infty)\times TQ$.
One such extension is obtained by composing $\phi$ with a sufficiently differentiable retraction~\cite[Ch.~6]{Lee2012} of $TQ$ onto $TA$ (such a retraction is guaranteed to exist locally by transversality of constraint surfaces).
}
$\phi|_{[0,\infty)\times\Sigma}$ is a continuous selection of the finite collection of $C^r$ functions 
$$\set{\phi_\omega^\eta : \omega\in\Omega, \eta\in H(\omega)}$$
on the open neighborhood 
$\e{F}\subset TQ$ 
containing $(t,(\q,\dot{\q}))$, 
i.e. 
$\phi|_{[0,\infty)\times\Sigma}$
is piecewise--$C^r$ at $(t,(\q,\dot{\q}))$.
See~\fig{hds} for an illustration the piecewise--differentiability of trajectory outcomes arising from a transition between contact mode sequences.
\end{proof}

\remna{relaxing Theorem hypotheses}{%
Since the class of piecewise--differentiable functions is closed under finite composition, conditions (i) and (ii) in the preceding Theorem can be readily relaxed to accommodate a finite number of constraint activation/deactivation times in the interval $(0,t)$.
Conditions (iii) and (iv) are more difficult to relax since there are systems wherein trajectories initialized arbitrarily close to an admissible trajectory fail to be admissible themselves.
As a familiar example, consider a 1 degree--of--freedom elastic impact oscillator~\cite[Ch.~2.4]{GuckenheimerHolmes1983} (i.e. a bouncing ball):
the stationary trajectory (initialized with $q,\dot{q} = 0$) is admissible for all time, but all nearby trajectories (initialized with $q \ne 0$ or $\dot{q}\ne 0$) exhibit the Zeno phenomenon.
We will discuss further possible extensions in~\sct{relax}.
}

Under the hypotheses of 
the preceding Theorem,
the continuous flow $\phi$ is piecewise--differentiable at a point $(t,(\q,\dot{\q}))\in [0,\infty)\times TA$,
that is,
near $(t,(\q,\dot{\q}))$ the graph of $\phi$ is covered by the graphs of a finite collection 
$\set{\phi_\omega^\eta : \omega\in\Omega, \eta\in H(\omega)}$
of differentiable functions (termed selection functions).
This implies in particular that there exists a continuous and piecewise--linear function 
\eqnn{
D\phi(t,(\q,\dot{\q})):
T_{(t,(\q,\dot{\q}))}\paren{[0,\infty)\times TA}
\into T_{\phi(t,(\q,\dot{\q}))}A 
}
(termed the Bouligand or B--derivative) that provides a first--order approximation for how trajectory outcomes vary with respect to initial conditions.
Formally, for all $(u,(v,\dot{v}))\in 
T_{(t,(\q,\dot{\q}))}\paren{[0,\infty)\times TA}%
$, 
the vector 
$D\phi(t,(\q,\dot{\q});u,(v,\dot{v}))\in\R^{2d}$ is the directional derivative of $\phi(t,(\q,\dot{\q}))$ in the $(u,(v,\dot{v}))$ direction:
\eqnn{
\label{eq:foa}
\lim_{\alpha\downarrow 0} \frac{1}{\alpha}\left[\paren{\phi(t+\alpha u,(\q+\alpha v, \dot{\q}+\alpha\dot{v})) 
		-\phi(t,(\q,\dot{\q}))} -\right. \\
	\left. D\phi(t,(\q,\dot{\q});u,(v,\dot{v})) \right] = 0.
}
Furthermore, this directional derivative is contained within the collection of directional derivatives of the selection functions.
Formally, for all $(u,(v,\dot{v}))\in 
T_{(t,(\q,\dot{\q}))}\paren{[0,\infty)\times TA}%
$,
\eqnn{
D\phi(t,(\q,\dot{\q});u,(v,\dot{v}))\in\quad\quad\quad\quad\quad\quad\quad\quad\quad\quad\quad\\
\set{D\phi_\omega^\eta(t,(\q,\dot{\q});u,(v,\dot{v})) : \omega\in\Omega,\eta\in H(\omega)}.
}
The selection functions are classically differentiable, whence their directional derivatives can be computed via matrix--vector multiplication between a classical (Jacobian/{\Frechet}) derivative matrix and the perturbation vector.
Formally, for all $(u,(v,\dot{v}))\in 
T_{(t,(\q,\dot{\q}))}\paren{[0,\infty)\times TA}%
$, 
$\omega\in\Omega$,
$\eta\in H(\omega)$,
\eqnn{
D\phi_\omega^\eta(t,(\q,\dot{\q});u,(v,\dot{v})) = D\phi_\omega^\eta(t,(\q,\dot{\q}))\mat{c}{u \\ v \\ \dot{v}},
}
where $D\phi_\omega^\eta(t,(\q,\dot{\q}))\in\R^{(2d)\times(1+2d)}$ 
is the classical derivative of the selection function $\phi_\omega^\eta$.
The matrix $D\phi_\omega^\eta(t,(\q,\dot{\q}))$ can be obtained by applying the (classical) chain rule to the definition of $\phi_\omega^\eta$ from~\eqref{eq:phiomegaeta}.

\section{Discussion}
\label{sec:disc}

We conclude by discussing possible routes (or obstacles) to extend our result,
and implications for assessing stability and controllability.

\subsection{Extending our result}
\subsubsection{Relaxing hypotheses}\label{sec:relax}
The hypotheses used to state {\thmdiffb} restrict either the systems or system trajectories under consideration; we will discuss the latter before addressing the former.

Trajectories we termed admissible exhibit neither \emph{grazing} nor \emph{Zeno} phenomena.
Since grazing generally entails constraint activation times that are not even Lipschitz continuous with respect to initial conditions, the flow is not piecewise--$C^r$ along grazing trajectories.
This fact has been shown by others~\cite[Ex.~2.7]{Di-BernardoBudd2008}, and 
is straightforward to see in an example.  
Indeed, consider the trajectory of a point mass moving vertically in a uniform gravitational field subject to a maximum height (i.e. ceiling) constraint.
The grazing trajectory is a parabola, whence the time--to--activation function involves a square root of the initial position.
Zeno trajectories, on the other hand, can exhibit differentiable trajectory outcomes 
following an accumulation of constraint activations (and, hence, deactivations); consider, for instance, the 
(stationary) outcome that follows the accumulation of impacts in a model of a 
bouncing ball~\cite[Ch.~2.4]{GuckenheimerHolmes1983}.
Thus we cannot at present draw any general conclusions regarding differentiability of 
the flow along Zeno trajectories, and speculate that it might be possible to recover 
piecewise--differentiability along such trajectories in the \emph{completion} of the mechanical 
system~\cite[Sec.~IV]{OrAmes2011} after establishing continuity with respect to 
initial conditions in the intrinsic state--space metric~\cite[Sec.~III]{BurdenGonzalezVasudevan2015tac}.

The criteria we impose on the mechanical system are more numerous.
{\asmporth} can be generalized to include bilateral constraints
by simply imposing orthogonality of the representation of the unilateral constraints in the intrinsic constraint manifold compatible with the bilateral constraints.
We restricted the configuration space to $Q = \R^d$ starting in~\sct{mdl} to simplify the exposition and lessen the (already substantial) notational overhead,
but the preceding results apply to more general configuration spaces:
if pairing $Q$ with the inertia tensor $M$ yields a Riemannian manifold with respect to which Assumptions 1, 2, and 4 are satisfied, 
then the preceding results remain valid.%
\footnote{Since the preceding results are not stated in coordinate--invariant terms, the differentiability results are formally applicable only after passing to coordinates.}

The dynamics in~\eqref{eq:dyn} vary with the contact mode $J\subset\set{1,\dots,n}$ due to intermittent activation of unilateral constraints $a_J(q) \ge 0$,
but the (so--called~\cite{Ballard2000}) effort map $f$ was not allowed to vary with the contact mode.
Contact--dependent effort can easily introduce nonexistence or nonuniqueness.
Indeed, this phenomenon was investigated thoroughly by Carath{\'{e}}odory and, later, Filippov~\cite[Ch.~1]{Filippov1988}.
For a specific example of the potential challenges in allowing contact-dependent forcing, note that the introduction of simple friction models into mechanical systems subject to unilateral constraints is known to produce pathologies including nonexistence and nonuniqueness of trajectories~\cite{Stewart2000}.
To generalize the preceding results to allow the above phenomena, one would need to 
provide conditions ensuring that trajectories 
(i) exist uniquely, 
(ii) depend continuously on initial conditions, 
and 
(iii) admit differentiable selection functions along trajectories of interest.

\subsubsection{Including control inputs}
\label{sec:control}
We focused on autonomous dynamics in~\eqref{eq:dyn}; however, parameterized control inputs can be incorporated through a standard state augmentation technique in such a way that Theorem~\ref{thm:diffb} implies trajectory outcomes depend piecewise--differentiably on initial states and input parameters, even as the contact mode sequence varies.

Specifically, suppose~\eqref{eq:dyn:cont} is replaced with
\begin{subequations}\label{eq:ctrl}
\begin{align}
 M(q)\ddot{q} & = \tilde{f}((q,\dot{q}), u) + c(q,\dot{q})\dot{q} + Da_J(q)^\tr \tilde{\lambda}_J((q,\dot{q}), u),\label{eq:ctrl:cont}\\
 \dot{q}^+ & = \td{\Delta}_J((q,\dot{q}^-),u)\dot{q}^-,\label{eq:ctrl:disc}
\end{align}
\end{subequations}
where $\tilde{f}\colon TQ \times U \into \R^d$ is an effort map that accepts a constant
input parameter $u\in U = \R^m$,
$\tilde{\lambda}_J\colon TQ \times U \into \R^{\abs{J}}$ is the reaction force that results from applying effort $\tilde{f}(q,\dot{q},u)$ in contact mode $J$,
and $\td{\Delta}_J\colon TQ\times U\into\R^d$ is a reset map that accepts input parameter $u$ as well.
We interpret the vector $u$ as parameterizing an open-- or closed--loop input to the system;
once initialized, $u$ remains constant.%
\footnote{A control policy represented using a universal function approximator such as an artificial neural network~\cite{LevineFinn2016, KumarTodorov2016} provides an example of a parameterized closed--loop input, while a control signal represented using a finite truncation of an expansion in a chosen basis~\cite{MombaurLongman2005, KuindersmaDeits2015} provides an example of a parameterized open-loop input.}
It is possible to generalize the proof of~{\thmdiffb} to provide conditions under which there exists a continuous flow $\td{\phi}:\td{\e{F}}\into TA$
for~\eqref{eq:ctrl} that is piecewise--differentiable with respect to initial conditions $(q,\dot{q})\in TA$ and input parameters $u\in U$ over an open subset $\td{\e{F}} \subset [0,\infty)\times TA\times U$ containing $\set{0}\times TA\times U$.

\subsection{Assessing (in)stability of periodic orbits}

\label{sec:periodic}
\begin{figure}[h!]
\centering
\includegraphics[width=.9\columnwidth]{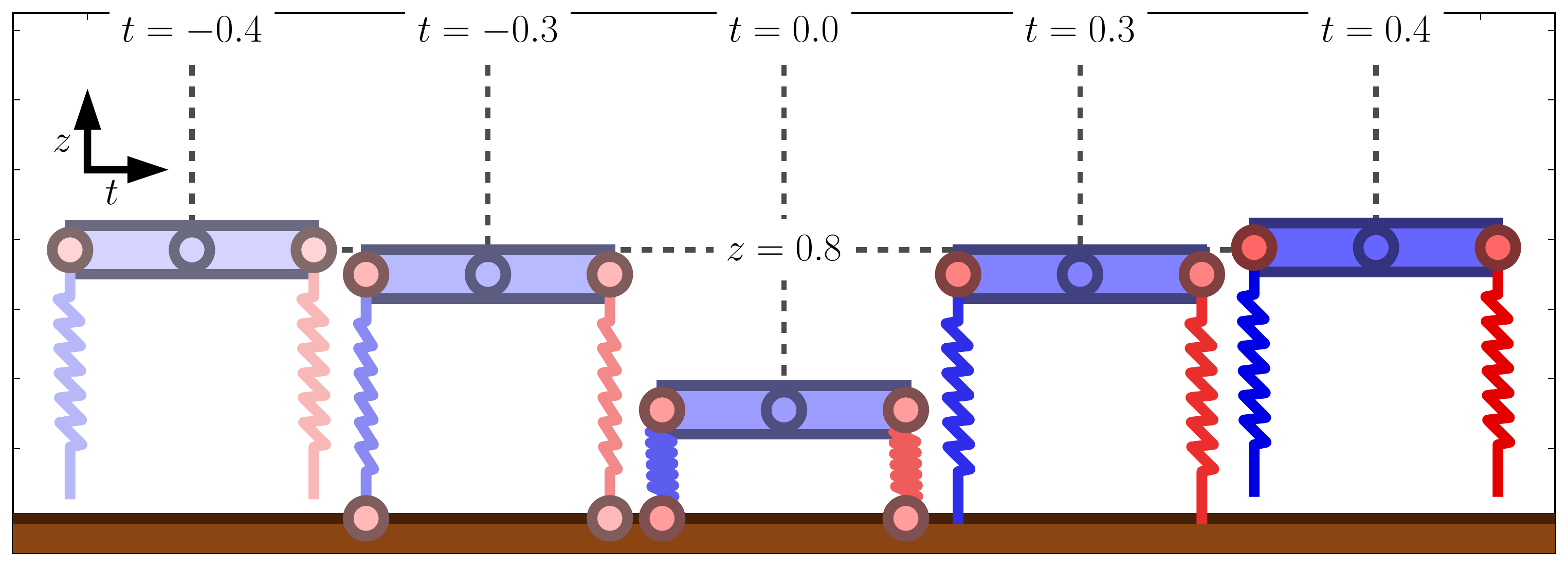}
\\
\resizebox{.9\columnwidth}{!}{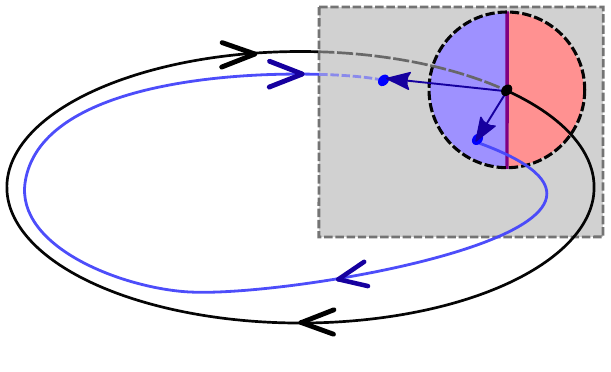}
\caption{
\label{fig:ex_poincare}
Illustration of a periodic orbit in the system depicted in Fig.~\ref{fig:ex}(\emph{right}) 
undergoing simultaneous activation (and, subsequently, simultaneous deactivation) of unilateral constraints.
(\emph{top})
Snapshots of the trajectory at apex ($t = \pm 0.4$), touchdown/liftoff ($t = \pm 0.3$), and nadir ($t = 0.0$).
(\emph{bottom})
Illustration of {\Poincare} map $P:V\into S$ over the apex section $S$,
defined in \sct{periodic};
trajectories initialized in the open subset $V\subset S$ return to $S$.
The set $V$ is partitioned into regions where selection functions for the piecewise--$C^r$ map $P$ are active:
initial conditions with $\theta>0$, where the left leg constraint activates before the right, are colored blue;
initial conditions with $\theta<0$ are colored red.
Along the trajectory generated by the fixed point $P(\rho,\dot{\rho}) = (\rho,\dot{\rho})$ (colored black),
simultaneous constraint activation is indicated by the trajectory passing through the intersection of the constraint surfaces for the right ($\set{a_r=0}$) and left ($\set{a_l=0}$) legs;
similarly for simultaneous deactivation through the intersection $\set{\lambda_l = 0}\cap\set{\lambda_r = 0}$.
A nearby trajectory initialized at $(p,\dot{p})\in V$ (colored blue) undergoes constraint activation and deactivation at distinct instants in time.
}
\end{figure}

In this section we consider the problem of assessing stability (or instability) of a periodic orbit in a mechanical system subject to unilateral constraints.
Suppose $(\rho,\dot{\rho})\in TA_{\emptyset}$ is an initial condition that lies on a \emph{periodic orbit}, i.e. there exists $T > 0$ so that $\phi(T,(\rho,\dot{\rho})) = (\rho,\dot{\rho})$ (and so that $\phi(t,(\rho,\dot{\rho})) \ne (\rho,\dot{\rho})$ for all $t\in(0,T)$).
If the trajectory $\phi^{(\rho,\dot{\rho})}$ undergoes constraint activations and deactivations at isolated instants in time, then prior work has shown that $\phi$ is $C^1$ at $(T,(\rho,\dot{\rho}))$, and the classical derivative $D\phi(T,(\rho,\dot{\rho}))$ can be used to assess stability of the periodic orbit~\cite{AizermanGantmacher1958}.
If instead the trajectory activates and/or deactivates some constraints simultaneously as in~\fig{ex_poincare}, then (so long as constraint activations/deactivations are admissible on and near $\phi^{(\rho,\dot{\rho})}$) 
the results of~\sct{diffa} ensure that $\phi$ is $PC^1$ at $(T,(\rho,\dot{\rho}))$ and the B--derivative $D\phi(T,(\rho,\dot{\rho}))$ is not generally given by a single linear map, whence classical tests for stability are not applicable.
In what follows we generalize the classical techniques to use this B--derivative to assess stability (or instability) of the periodic orbit $\phi^{(\rho,\dot{\rho})}$.

We start by constructing a {\Poincare} map for the periodic orbit $\phi^{(\rho,\dot{\rho})}$.
Let $S \subset TQ$ be a {\Poincare} section for $\phi^{(\rho,\dot{\rho})}$ at $(\rho,\dot{\rho})$, 
i.e. a $C^r$ embedded codimension--1 submanifold containing $(\rho,\dot{\rho})$ that is transverse to the vector field in~\eqref{eq:dyn:cont}. 
For a concrete example we refer to the model in~\fig{ex_poincare}
where $S$ is a {\Poincare} section about an apex height of $z~=~0.8$ and
$\rho$ is a position vector with body height $z=0.8$, rotation $\theta = 0$, and the legs oriented perpendicular to the body orientation. 
Given zero initial velocity, the time period is $T=0.8$.

Since $\phi$ is continuous by~{\lemcont},
there exists a \emph{first--return time} $\tau \colon V \into (0, \infty)$ 
defined over an open neighborhood $V\subset S$ containing $(\rho, \dot{\rho})$ 
such that $\phi(\tau(q, \dot{q}), (q, \dot{q})) \in S$ for all $(q, \dot{q}) \in V$ and $\tau(\rho,\dot{\rho}) = T$;
we let $P:V\into S$ be the {\Poincare} (or \emph{first--return}) map defined by
\eqnn{
\forall (q, \dot{q}) \in V \colon P(q, \dot{q}) = \phi(\tau(q, \dot{q}), (q, \dot{q})) \in S.
}
As an illustration, $(p, \dot{p}) \in V$ in \figref{fig:ex_poincare}(\emph{bottom}) generates a trajectory initialized near $(\rho,\dot{\rho})$ that undergoes constraint activations and deactivations at distinct instants in time,
activating the left leg constraint before activating the right leg constraint,
then deactivating both constraints in the same order.
Since $\phi$ is $PC^r$ and $S$ is a $C^r$ manifold we conclude that $\tau$ is $PC^r$~{\cite[Thm.~10]{BurdenSastry2016siads}}, whence $P$ is $PC^r$.
This implies in particular that its B--derivative $DP(\rho,\dot{\rho})$ provides a continuous and piecewise--linear first--order approximation for $P$.
To assess exponential stability of $\phi^{(\rho,\dot{\rho})}$, it suffices to determine conditions under which the piecewise--linear map $DP(\rho,\dot{\rho})$ is exponentially contractive or expansive. 
This task is nontrivial since, as is well--known~\cite[Ex.~2.1]{Branicky1998}, a piecewise--linear system constructed from stable subsystems may be unstable; 
similarly, a system constructed from unstable subsystems may be stable.
We refer to~\cite[Sec.~II-A]{LinAntsaklis2009} for a thorough review of state--of--the--art methods for assessing stability of piecewise--linear systems, and provide some example tests below.

\begin{figure}[h]
\centering
{
\def\svgwidth{0.5\columnwidth}
\resizebox{.45\columnwidth}{!}{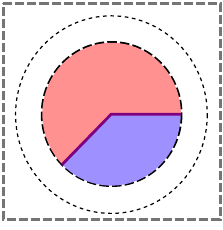}
\quad
\def\svgwidth{0.5\columnwidth}
\resizebox{.45\columnwidth}{!}{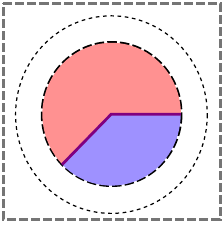}
\caption{
\label{fig:stability}
A \Poincare section $S$ (outer grey box) and the neighborhood $V$ (inner circle) containing $(\rho,\dot{\rho})$ over which the piecewise--differentiable {\Poincare} (or first--return) map $P:V\into S$ is defined.
In this example, the color corresponds to one of the two selection functions and the shaded region
of $V$ denoting the active selection function for a perturbation in the given direction.
$B$ is the unit ball.
The dotted ellipses are the unit ball transformed by the respective selection functions with the arrow indicating the principle axes.
(\emph{left})
The \Poincare map $P$ corresponds to an unstable system. One of the eigenvectors $\nu_1$ maps a unit vector outside the unit ball. Additionally, the selection function for which $\nu_1$ is an eigenvector for is active for any perturbation lying along $\nu_1$.
(\emph{right})
The given instability check is not able to determine if the \Poincare map $P$ is unstable due to the 
lack of active eigenvectors mapping a unit vector outside the unit ball.
The active selection function for a perturbation along $\nu_2$ is not the selection function for which $\nu_2$
is an eigenvector for.
}
}
\end{figure}

Since $P$ is $PC^r$,
there exists a finite collection $\set{P_\w}_{\w\in \Omega}$ of $C^r$ selection functions for $P$,
and we assume the neighborhood $V$ was chosen sufficiently small that $P_\omega\colon V\into S$ for each $\omega\in\Omega$.
Let $R_{\w} \subset V$ denote the region where the selection function $p_\w$ is \emph{active} (i.e. where $P|_{R_{\omega}} = P_\omega|_{R_{\omega}}$).
The first order approximation for $P_\w$ is given by the classical (Jacobian/\Frechet) derivative $DP_\w\colon TV \into TS$, which can be calculated using the (classical) chain rule.
If there is a norm $\norm{\cdot}:\R^{2d-1}\into\R$ 
with respect to which $DP_\w(\rho, \dot{\rho})$ is a contraction for all $\omega\in\Omega$ 
(i.e. for all $\omega\in\Omega$ the induced norm $\norm{DP_\w(\rho, \dot{\rho})} < 1$),
then the periodic orbit $\phi^{(\rho,\dot{\rho})}$ is exponentially stable~\cite[Prop.~15]{BurdenSastry2016siads}.
(Note that it does not suffice to find a different norm $\|\cdot\|_\omega$ for each $\omega\in\Omega$ with respect to which $DP_\omega(\rho,\dot{\rho})$ is a contraction.~\cite[Ex.~2.1]{Branicky1998}).
If instead for some $\omega\in\Omega$ there exists an eigenvector $\nu$ for $DP_\w(\rho, \dot{\rho})$ with eigenvalue $\lambda$ such that $\abs{\lambda} > 1$ and both
$\nu$ and $DP_i(\rho,\dot{\rho})\nu \in R_{i}$,
then $(\rho,\dot{\rho})$ is exponentially unstable;
this instability test is illustrated in~\figref{fig:stability}.

\subsection{Assessing controllability}
In this section we consider the problem of assessing (\emph{small--time}, \emph{local}~\cite{Sussmann1987}) controllability along a trajectory in a mechanical system subject to unilateral constraints.
The local control problem has been solved quite satisfactorily along trajectories in such systems that undergo constraint activation and deactivation at distinct instants in time for cases where the control input influences the discrete--time~\cite{LongMurphey2011} or continuous--time~\cite{RijnenSaccon2015} portions of~\eqref{eq:dyn}. 
We concern ourselves here with the controlled dynamics in~\eqref{eq:ctrl}, and focus our attention on trajectories that activate and/or deactivate multiple constraints simultaneously since (to the best of our knowledge) this case has not previously been addressed in the literature.

Toward that end, let $\td{\phi}:\td{\e{F}}\into TA$ be the flow of~\eqref{eq:ctrl} (a mechanical system subject to unilateral constraints with input parameter $u\in U = \R^m$), 
and let $\td{\phi}^{((\rho,\dot{\rho}),\mu)}$ be a trajectory initialized at $(\rho,\dot{\rho})\in TA$ with input parameter $\mu\in U$.
If $\td{\phi}$ were $C^1$ at $(t,(\rho,\dot{\rho}),\mu)\in\td{\e{F}}$, then (small--time) local controllability about $\td{\phi}^{((\rho,\dot{\rho}),\mu)}$ could be determined using an invertibility condition on the (Jacobian) matrix $D\td{\phi}(t,(\rho,\dot{\rho}),\mu)$.
Indeed, a straightforward application of the Implicit Function Theorem~\cite[Thm.~C.40]{Lee2012} 
shows that if the subblock $D_U\td{\phi}(t,(\rho,\dot{\rho}),\mu)$, 
which transforms first--order variations in the input parameter $u$ into the resulting first--order variations in the state $(q,\dot{q})$ at time $t$, is invertible,
then~\eqref{eq:ctrl} is (small--time) locally controllable along $\td{\phi}^{((\rho,\dot{\rho}),\mu)}$~\cite[Thm.~8]{LevinNarendra1993}.%
\footnote{It will be useful in what follows to note that this invertibility condition is equivalent to the existence of a linear homeomorphism relating variations in (an appropriately--chosen subspace of) input parameters to variations in system states.}

In contrast to the preceding discussion, suppose now that $\td{\phi}^{((\rho,\dot{\rho}),\mu)}$ undergoes simultaneous constraint activations in the time interval $(0,t)\subset[0,\infty)$. 
In this case $\td{\phi}$ will not be $C^1$ at $(t,(\rho,\dot{\rho}),\mu)$, so the classical test for controllability is not applicable. 
If all constraint activations and deactivations are admissible for $\td{\phi}^{((\rho,\dot{\rho}),\mu)}$ and nearby trajectories,
then {\thmdiffb} implies that $\td{\phi}$ is $PC^r$ at $(t,(\rho,\dot{\rho}),\mu)$ and hence possesses a B--derivative $D\td{\phi}(t,(\rho,\dot{\rho}),\mu)$, that is, a continuous and piecewise--linear first--order approximation.
By analogy with the classical test~\cite[Thm.~8]{LevinNarendra1993}, a variant of the Implicit Function Theorem applicable to $PC^r$ functions~\cite[Thm.~4.2.3]{Scholtes2012} can be used to derive a sufficient condition for small--time local controllability along $\td{\phi}^{((\rho,\dot{\rho}),\mu)}$:
if the piecewise--linear function that transforms first--order variations in (an appropriately--chosen subspace of) input parameters $u$ into the resulting first--order variations in the state $(q,\dot{q})$ at time $t$ is a (piecewise--linear) homeomorphism, 
then~\eqref{eq:ctrl} is (small--time) locally controllable along $\td{\phi}^{((\rho,\dot{\rho}),\mu)}$. 

\iftoggle{acm}{
\renewcommand*{\bibfont}{\small}
}
{
\pagebreak
\hypersetup{linkcolor=blue}
}
\printbibliography

\end{document}